\documentclass[11pt]{article}
\usepackage[latin1]{inputenc}
\usepackage{amsmath,amssymb}
\usepackage{color}
\usepackage{multicol}
\usepackage{ulem,stmaryrd}
\usepackage{latexsym,epsfig}
\usepackage{amssymb,amsmath,amsfonts,amscd}
\usepackage{exscale}
\usepackage{ifthen}
\usepackage{longtable}
\usepackage{subfigure}
\usepackage{cite}
\usepackage{lscape}
\usepackage{soul}
\usepackage{ulem}
\usepackage{cancel}

\newtheorem{theo}{\textbf{Theorem}\ }
[section]
\newtheorem{lemma}[theo]{\textbf{Lemma}\ }
\newtheorem{coro}[theo]{Corollary\ }

\newtheorem{prop}[theo]{\textbf{Proposition}\ }
\newtheorem{properties}[theo]{\textbf{Properties}\ }

\numberwithin{equation}{section}   

\begin{document}

	\begin{center}
		{\Large \bf The survival probability }
		
			\vspace{0.3cm}
			 
		{\Large \bf of a weakly subcritical multitype branching process}
		
			\vspace{0.3cm}
			
		{\Large \bf  in iid random environment.} 
		
		 \vspace{1cm}
		
		{ M. Peign\'e} $^($\footnote{
			Institut Denis Poisson UMR 7013,  Universit\'e de Tours, Universit\'e d'Orl\'eans, CNRS  France.  email: peigne@univ-tours.fr,}$^) \ \&$ 
		 		{ C. Pham} $^($\footnote{
			CERADE ESAIP, 18 rue du 8 mai 1945 - CS 80022 - 49180 St-Barth\'elemy d'Anjou ; email: dpham@esaip.org}$^), ^($\footnote{
			LAREMA UMR CNRS 6093, Angers, France. }$^)$

		\today
		
	\end{center}
	
	\noindent {\bf  \small Abstract. } 
	 We study the asymptotic behavior of the  probability of  non extinction of a weakly subcritical multitype branching  process in iid random environment. 
        Under suitable assumptions, the survival  probability is of order of  $\rho^n n ^{-3/2}$ for some $\rho \in (0, 1)$ to specify.

	\vspace{0.2cm}

	\vspace{0.5cm}

	
	
	\section{Introduction and main results}
%
 	\subsection{Historical context}
 		Branching processes in random environment (BPRE) is an important subject of the theory of branching processes. This model was introduced at the beginning of the 1960's; it takes into account the random fluctuations of the reproduction laws over time.  Various properties of this and more general models of BPRE  have  been analyzed these last decades, through a large number of varied articles.  There exists a relatively complete description of the basic properties of many models of BPRE either under the annealed approach or the  quenched one. In particular, the behavior of the single type BPRE is mainly determined be the properties of the so-called ``associated random walk'' constructed by the logarithms of the expected population sizes of particles of the respective generations; this random walk   divides in a natural way the set of all single type BPRE into the classes of {\it supercritical}, {\it critical} and {\it subcritical}  processes (see [6] for more detail).  Notice that a precise estimate for single type branching process in finite state space markovian environment does exist \cite{GLL2019}; the fact that the underlying Markov chain is finite is essential in this study, this hypothesis is not relevant in the case of product of random matrices.

		Analogue  statements known for the single type case do exist for the multitype BPRE with finitely types $p\geq 2$. The role of the associated random walk is played in this case by the logarithms of the norms of products of the mean random matrices associated to the environment. Recent useful results on fluctuations of  ordinary random walks on the real line have been extended  in terms of $p \times p$ random matrices. In particular, the description of  the asymptotic behavior of the survival probability for the critical multitype BPRE under general conditions is done in \cite{LPP2018} and  \cite{DV2018} where it is proved, under natural and quite general assumptions, that the probability of survival up to time $n$  is of  order $ n^{1/2}$ for a  multitype branching process evolving in iid random environment.  The case of supercritical multitype branching processes is  studied in \cite{GLP2022}, where a Kesten-Stigum type theorem is established.  At last,  as in the single type case,   subcritical BPRE are divided in three categories: {\it strongly}, {\it intermediately} and {\it weakly} subcritical. The first two subcases  are the subject of recent work:  in the strongly subcritical subcase, the annealed survival probability at time $n$   is equivalent to $\rho^n$ \cite{VW} while in the intermediately subcritical subcase it is  of  order   $\rho^n n ^{-1/2}$ \cite {DV2020}, where $\rho \in(0, 1)$ is a constant defined in terms of the Lyapunov exponent for products of the mean value matrices of the laws of reproduction of particles.  		
		
		In the present paper, we establish a rough asymptotic estimate of the annealed survival probability at time $n$ in the weakly subcritical case. This case was  resistant since there was no local limit theorem for the norm of products of random matrices conditioned to remain greater than 1 until time $n$.  Such a local limit theorem is  obtained in  \cite{PP2023}, its proof can be adapted to the present context  to show that the probability of survival at time $n$ is of order 
		$\rho^n n ^{-3/2}$ where $\rho \in (0, 1)$. 
 \subsection{Notations and assumptions}

We fix an integer $p\geq 2$ and  denote $\mathbb R^p$ (resp.,  $\mathbb N^p$) the set of $p$-dimensional column vectors with real (resp.,    non negative integers) coordinates; for any column vector $x = (x_i)_{1\leq i\leq p} \in \mathbb R^p$, we denote $\tilde x$ the row vector $\tilde x:= (x_1, \ldots, x_p)$. Let  ${\bf 1}$ (resp.,  ${\bf 0}$)  be the column vector of $\mathbb R^p$ whose all coordinates equal  $1$ (resp., $0$).
We fix  a basis $ \left\{{e}_i, 1 \le i \le p\right\}$  in $\mathbb{R}^p$ and denote $\vert \cdot \vert $ the corresponding $ \mathbb L^1$ norm. 

The multitype Galton-Watson process is a temporally homogeneous vector Markov process  $(Z_n)_{n \geq 0}$  whose states are column vectors in $\mathbb{N}^p$.   For any $1\leq i\leq p$,   the $i$-th component $Z_n (i)$ of $Z_n$ may be interpreted   as the number of objects of type $i$ in the $n$-th generation.

A multivariate  probability generating function  $ 
 (f^{(i)} )_{1\leq i \leq p}$, denote by $f$,  is a function from $(\mathbb R^{+})^p $ to $\mathbb R^+$ defined by  
$$
f^{(i)} (s) = \sum_{ \alpha  \in {\mathbb{N}^p}} {{p^{(i)}   }(\alpha)}s^{ \alpha}, 
$$
for any $s=(s_i)_{1\leq i\leq p}\in (\mathbb R^+)^p$, where
\begin{enumerate}
	\item [(i)] $\alpha =(\alpha_i)_i  \in \mathbb{N}^p$ and $s^{  \alpha} = s_1^{{\alpha _1}} \ldots s_p^{{\alpha _p}}$;
	\item [(ii)]$p^{(i)}  (\alpha)={p^{(i)}  } ({{\alpha _1}, \ldots,{\alpha _p}})$ is the  probability that an object of type $i$   has $\alpha_1$ children of type $1,\ldots,\alpha_p$ children of type $p$.
\end{enumerate}

From now on, the multivariate generating function $f$ is a random variable defined on  a  probability space $ ({\Omega,\mathcal{F},\mathbb{P}})$; let {\bf f}= $(f_n)_{n \geq 0}$ be a sequence of   iid random copies of $f$, called the {\it (random) environment}. The Galton-Watson process with $p$ types of particles in  the random environment {\bf f} describes the evolution of a particle population ${Z_n} =  (Z_n(i))_{1\leq i \leq p}$  for $n \geq 0$.
 
For any $ s=(s_i)_{1\leq i\leq p}, 0 \le s_i \le 1$,
	\begin{eqnarray*} \label{eqn1.1}
		\mathbb{E} \left( s^{Z_n} | Z_0, \ldots, Z_{n-1}, f_0, \ldots, f_{n-1} \right) = f_{n-1}(s) ^{Z_{n-1}}  
	\end{eqnarray*} 
	which yields  
\[
	  \mathbb{E} \left(s^{Z_n}| Z_0= {\tilde e}_i, f_0, \ldots, f_{n-1} \right) =
		  f_0^{(i)} (f_1 (\ldots  f_{n-1}(s)\ldots)).
\]
The  probability of  non extinction $q_n^{(i)}$ at generation $n$ given the environment ${\bf f}$ when the ancestor is of type $i$  is   
	\begin{eqnarray*} \label{eqn2.11}
		q_n^{(i)} &:=& \mathbb{P} (\vert Z_n \vert >0\mid f_0^{(i)}, f_1,\ldots, f_{n-1}) \notag \\
		&=& 1-f_0^{(i)} (f_1 (\ldots  f_{n-1}({\bf 0})\ldots)) = \tilde e _i  ({\bf{1}}- f_0 (f_1 (\ldots  f_{n-1}({\bf 0})\ldots))),
	\end{eqnarray*}
	so that
	\begin{eqnarray*}
	\mathbb E [q_n^{(i)}] =   \mathbb{P}(Z_n \ne \tilde {\bf 0}| Z_0= \tilde e _i)= \mathbb E[\tilde e _i  ({\bf{1}}- f_0 (f_1 (\ldots  f_{n-1}({\bf 0})\ldots)))].
	\end{eqnarray*}
The asymptotic behavior of the quantity above is controlled by the mean of the offspring distributions.   
We assume that the offspring distributions have finite first and second moments;  the generating function  $f = (f^{(i)})_{ 1\leq i\leq p}$, is thus $C^2$-functions on $[0, 1]^p$ and we introduce:
\begin{enumerate}
\item [(i)]
 the  mean matrix $M = (M (i, j))_{i, j}=   \left({\displaystyle \frac{{\partial f^{(i)}  }}{{\partial {s_j}}}} ({\bf 1})\right)_{i,j}$; 

\item  [(ii)] the  Hessian matrices $B ^{(i)}=    \left(\displaystyle{\partial^2f^{(i)}  \over \partial s_k\partial s_\ell}(\mathbf 1) \right)_{k, \ell},\quad \text{for}\ i=1, \ldots, p$. \end{enumerate}

We denote by $(M_n)_{n \geq 0}$ (resp. $(B_n^{(i)})_{n \geq 0}, i=1, \ldots, p)$ the sequence of iid random mean matrices (resp. Hessian matrices) associated with the sequence $(f_n)_{n \geq 0}$. These matrices belong to the semigroup $\mathcal S$ of $p\times p$ matrices with positive entries; we endow $\mathcal S$ with the $\mathbb L^1$-norm.

We  assume that the distribution $\mu$ of $M$ satisfies the following assumptions.

	\noindent {\bf P1} {\it  \underline{Moment assumption}:   $\mathbb E   
	(\vert M\vert) <+\infty.$
	}

	\noindent {\bf P2}  {\it  \underline{Irreducibility  assumption}: The support of  $\mu$  acts strongly irreducibly on the semigroup of matrices with  non negative entries, i.e. there exists no affine subspaces $A$ of $\mathbb R^p$ such that $A \cap (\mathbb R^+)^p$ is  non empty, bounded and invariant under the action of all elements of the support of $\mu$.}

	\vspace{2mm}
	
	\noindent {\bf P3} 
  {\it There exists $B>1$ such that for  $\mu$-almost all $M$    and any $1\leq i, j, k, \ell\leq p,$ } 
\begin{equation}\label{sdelta}
{1\over B} M(k, \ell) \leq M(i, j) \leq B \    M(k, \ell).
\end{equation} 
 From now on, we denote by $\mathcal S_{B}$ the subset of $\mathcal S$  of $p\times p$-matrices $M$ satisfying condition (\ref{sdelta}).	
	   		
\noindent {\bf P4} {\it  \underline{Subcriticality}: The upper Lyapunov exponent $\displaystyle \gamma_\mu:= \lim_{n \to +\infty} {\mathbb E[\ln \vert M_0\cdots M_{n-1}\vert ]\over n}$ is negative.}

	\noindent {\bf P5} {\it There exists  $\varepsilon>0$  such that 
	$
	\mu \{M \in \mathcal S: \forall x \in (\mathbb R^+)^p \ {\rm such that } \ \vert x\vert =1, \ln \vert  \tilde x M \vert   \geq {\tiny \varepsilon} \} >0
	 $
	 and 
	  $
	 \mu \{M \in \mathcal S: \forall x \in (\mathbb R^+)^p \ {\rm such that } \ \vert x\vert =1, \ln \vert   M x \vert   \geq \varepsilon \} >0
	 $.} 
	
	\vspace{2mm}
	   As it is usual in studying local probabilities, one has to distinguish between ``lattice'' and ``non lattice'' distributions. The ``non lattice'' assumption   ensures that the  process $(\ln \vert \tilde x M_{0, n-1}\vert)_{n \geq 1}$ does not live in the translation of a proper subgroup of $\mathbb  R$; in the contrary case, when $\mu$ is lattice,  a phenomenon of cyclic classes appears  which involves some complications which are not interesting in our context.   In  section \ref{productmatrices}, we  give a  precise definition of this notion  in the context of products of random matrices.

 \noindent {\bf P6} {\it \underline{ Non lattice assumption}:  The measure $\mu$ is  non lattice.}
		
		\vspace{2mm} We also introduce the following hypotheses concerning  the environment ${\bf f}$.		
		\vspace{2mm}
		
		\noindent {\bf P7} {\it There exist $\varepsilon_0, K_0>0 $  such that for any $i, j \in \{1, \ldots, p\}$ and any $n \geq 0$, }
\begin{align*}
&(a) \ \mathbb P(Z_{n+1}(i) \geq 2\mid Z_{n}= \tilde e_j)\geq \varepsilon_0, \\
&  (b) \ \mathbb P(\vert Z_{n+1}\vert =0\mid Z_{n}= \tilde e_j)\geq \varepsilon_0,\\
&  (c) \ \mathbb E[\vert Z_{n+1}\vert^2 \mid Z_{n}= \tilde e_j]\leq K_0. 
\end{align*}
		\vspace{2mm}
		

Let us first explain some consequences of hypotheses {\bf P1} and 	{\bf P4} and the way three different subcases appear in the subcritical case.
 Using the standard subadditivity argument, one infers that,   under   {\bf P1},  for any $\theta \in [0, 1]$ the limit
\[\lambda(\theta):=\lim_{n \to +\infty} \mathbb E \left[\left\vert M_{n-1}\cdots M_0\right\vert^\theta\right]^{1/n}<+\infty\] is well defined. The function $\Lambda$ is the analogue of the moment generating function for the associated random walk in the case of branching processes in random environment with single type of particles.
The function $\theta \mapsto \mathbb E[\vert M \vert ^\theta]$   is  at least twice continuously differentiable on $[0, 1[$.
We set: for   $\theta \in [0, 1[$,
\[\Lambda(\theta):= \ln \lambda(\theta).\]
This function is  also  twice continuously differentiable  and strictly convex on $[0, 1[ $ (see \cite{BDGM}, Proposition 3.1).
 There are three cases to consider.

1. \underline{$ \Lambda'(1) <0$} (hence $\Lambda(1) <0$, by convexity). This case  corresponds to the {\it strongly subcritical } case for single type branching processes, it is studied in \cite{VW} when $p\geq 2$. Under suitable conditions, it holds
\[
\mathbb E [q_n^{(i)}]=\mathbb P(\vert Z_n\vert >0\mid Z_0= \tilde e_i)\sim c^i\  \lambda(1)^n, \quad i=1, \ldots, p
\]
with $c^i>0$ and $\lambda(1) \in ]0, 1[$.

2. \underline{$\Lambda'(1) =0$} (hence $\Lambda(1) <0$). It corresponds to the  {\it intermediately subcritical } case  for single type branching processes, it is studied in \cite{DV2020} when $p \geq 2$. Under suitable conditions, it holds the following rough estimate
\[
\mathbb E [q_n^{(i)}]=\mathbb P(\vert Z_n\vert >0\mid Z_0= \tilde e_i)\asymp c^i \ {\lambda(1)^n\over \sqrt{n}}, \quad i=1, \ldots, p
\]
with $c^i>0$ and $\lambda(1) \in ]0, 1[$.

3. \underline{$ \Lambda'(1) >0$}. This is the {\it weakly subcritical } case when $p=1$ and this is the aim of the present article in the case when $p\geq 2$, we only obtain a rough estimate.

\vspace{3mm}

We denote by $\theta_\star$ the unique value of $\theta$ in $]0, 1[$   such that    $\Lambda'(\theta_\star)= 0$. We set $\rho_\star = \lambda (\theta_\star)$; it is obvious that $\rho_\star \in ]0, 1[$.

\vspace{2mm}

\noindent  
{\bf Notation.}
{\it  
	Let $c >0$ and $\phi, \psi$ be two functions of some variable $u$; we shall write
	$\phi\stackrel{c}{\preceq}\psi $  (or simply $\phi \preceq \psi $) when $\phi(u) \leq c \psi(u)$ for any value of $u$.
	The notation
	$\phi \stackrel{c}{\asymp}\psi $ (or simply $\phi \asymp \psi $) means  $\phi \stackrel{c}{\preceq}\psi \stackrel{c}{\preceq}\phi.$}

\subsection{Main statement}

Now, we may state the main result of this article which concerns a rough estimate of the probability of extinction in the weakly subcritical case.



 \begin{theo}\label{theoprincipal} Assume that conditions {\bf P1}-{\bf P7}. 
Assume that $\Lambda'(1)>0 $  and let $\theta_\star$ be  the unique value of $\theta \in ]0, 1[$ such that  $\Lambda'(\theta_\star)= 0$.
 
  Then, there exists positive constants ${\bf c}$ and  ${\bf C}$ such that for $1\leq i \leq p$,
 \[
 {\bf c}  \ { \rho_\star^n\over n^{3/2}}\leq \mathbb P(\vert Z_n\vert >0\mid Z_0={\tilde e}_i)\leq  {\bf C} \  { \rho_\star ^n\over n^{3/2}} 
 \]
 with $\rho_\star = \lambda(\theta_\star)\in ]0, 1[$.
 \end{theo}
 
	\section{On products of positive random matrices}\label{productmatrices}

	The random variables $M_{n} $ and  $B^{(i)}_{n} $ are iid  with values in the semigroup $\mathcal S$   of $ p \times p$ matrices with    positive  coefficients.  
	We set, for any $0\leq k\leq n, $
	\[M_{k, n}= M_k \cdots M_n \quad {\rm and } \quad M_{n, k} = M_n \cdots M_k.
	\]
%

  In order to control the asymptotic behavior of the matrices $M_{n, 0}$ (resp. $ M_{0, n}$), we  study their action   on  the cone   of column (resp. row) vectors    with positive entries. 

Let $\mathcal C$ (resp. $\tilde{ \mathcal C}$) be the set of column vector (resp. row vectors) in $\mathbb R^p$ with positive entries. For any $x \in \mathcal C$, we denote by $\tilde x$ the corresponding  row vector.   We also  set   $ \mathbb X := \{ x \in \mathcal C,  |x| =1 \}$ 
	and $\tilde{ \mathbb X }= \{\tilde x  \mid x \in \mathbb X\}$.
	 We consider the following (left and right) actions of $\mathcal S$:
	\begin{itemize}
		\item the  linear action of $\mathcal S$ on $\mathcal C$ (resp. $\tilde {\mathcal C}$) defined by $(M, x) \mapsto M x$ 	(resp.  $(M, \tilde x) \mapsto \tilde x M)$ for any  $M \in \mathcal S $ and $x \in \mathcal C$, 
		\item the  projective action of $\mathcal S$ on $\mathbb X$  (resp. $\tilde{\mathbb X}$) defined by $(M, x) \mapsto M \cdot x := \displaystyle \frac{M x}{\vert Mx\vert}$ 
		
		\noindent $\left({\rm resp.} \    (M,  \tilde x) \mapsto  \tilde x \cdot M = \displaystyle\frac{\tilde  x M}{\vert \tilde x M\vert}\right)$  for any  $M \in \mathcal S $ and $x \in \mathbb X$.
	\end{itemize}

\noindent For any $M   \in \mathcal S$,    set   
 $\displaystyle v( M) := \min_{1\leq j\leq d}\Bigl(\sum_{i=1}^p M(i, j)\Bigr).
 $
 Then,  for any $  x \in \mathcal C$,
 \begin{equation*}\label{controlnormgx}
 0< v(M)\  \vert x\vert \leq \vert Mx\vert \leq \vert  M\vert \ \vert x\vert.
 \end{equation*}
 Similarly, noticing $^t M$ the transpose  of the matrix $M$,  it holds 
 \begin{equation*}\label{controlnormgx}
 0< v(^tM)\  \vert x\vert \leq \vert \tilde x M\vert \leq \vert  M\vert \ \vert   x\vert.
 \end{equation*}
Consequently, hypothesis {\bf P5} holds when $\mu\{M\mid v(M) >1\}>0$ and $\mu\{M\mid v(^tM) >1\}>0$.

In the context of  branching processes, we   focus on the  right action of $\mathcal S$ on $\tilde{\mathbb X}$. Therefore, we naturally endow  $\tilde{\mathbb X}$  with a distance $d$ which is a variant of the Hilbert metric: this distance  is bounded on $\tilde{\mathbb X}$ and any element $M$ in $\mathcal S$ acts on the metric space $(\tilde{\mathbb X}, d)$ as a contraction.  
More precisely: for any $  \tilde x=(x_i)_{1\leq i\leq p},   \tilde y  = (y_i)_{1\leq i\leq p} \in \tilde{\mathbb X}$, we write
$$
m( \tilde x,  \tilde y)  = \min \Bigl\{ {x_i\over y_i} \Big\vert i = 1, \ldots, p  \ {\rm such \ that} \ y_i>0\Bigr\}
$$
and we set 
$
d(  \tilde x,   \tilde y) := \varphi\Bigl(m( \tilde x,  \tilde y)m(  \tilde y,   \tilde x)\Bigr),
$
where  $\displaystyle \varphi(s) := {1-s \over 1+s}$ for $0 \leq s\leq 1$.
For $M\in \mathcal S$, we set 
$$
c(M):= \sup \{d(\tilde x \cdot M, \tilde y \cdot M)\mid   \tilde x,   \tilde y \in  \tilde{\mathbb X}\}.
$$  By \cite{H97}, the function $d$ is a distance on $\tilde{\mathbb X}$ which satisfies the following properties.
\begin{properties} \label{propertiesHennion}
\begin{enumerate}
\item $\sup\{d(  \tilde x,  \tilde y)\mid  \tilde x,  \tilde y \in \tilde{\mathbb X}\}=1$.
\item For any $M\in \mathcal S, $
$$
c(M)= \max_{i,j,k,l \in \{1, \ldots, p\}} {\vert M(i, j)M(k, l)-M(i, l)M(k, j) \vert \over  M(i, j)M(k, l)+M(i, l)M(k, j)}.
$$
In particular, there exists $c_{B} \in [0, 1)$ such that $c(M) \leq c_{B} <1$ for any $M\in \mathcal S_{B}$.
\item $d(\tilde x\cdot M,  \tilde y\cdot M) \leq c(M) d(  \tilde x,   \tilde y) \leq c(M)$ for any $ \tilde x,   \tilde y \in \tilde{\mathbb X}$ and $ M \in \mathcal S_{B}$.
\item
$c(MM') \leq c(M) c(M')$ for any $M, M'   \in \mathcal S_{B}$.
\end{enumerate}
 \end{properties}

The following lemma plays a crucial role in controlling the behavior of the norm of the product of random matrices $M_0, M_1,  \ldots $ (see Lemma 2.2 in \cite{BPP}). We denote by $T_B$ the closed semigroup  generated by  $\mathcal S_B$.
 \begin{lemma}\label{keylemma}
 Under hypothesis {\bf P3}, for any $  M \in T_{B}$  and $1 \le i,j,k, \ell \le p,$
	\begin{eqnarray*} \label{eqn21}
	M(i,j) \stackrel{B^2}{\asymp} M(k,\ell).
	\end{eqnarray*}
		In particular, there exists $ \delta >1$ such that for any $ M, N \in T_B$ and   $  x,  y \in \mathbb X, $
		\begin{enumerate}
			\item  $  \vert { M x} \vert  \stackrel{\delta}{\asymp} \vert M \vert  \,\, \mbox{and} \,\, \vert \tilde y  M\vert  \stackrel{\delta}{\asymp} \vert  M\vert $, 
				\item $\vert \widetilde yMx\vert \,\, \stackrel{\delta}{\asymp} \,\vert M\vert $,
					\item $ \vert  M\vert \vert  N\vert \stackrel{\delta}{\preceq} \vert MN\vert  \leq \vert  M\vert \vert  N\vert $.
		\end{enumerate}
	\end{lemma}

 On the product space $\tilde {\mathbb X}  \times \mathcal S   $,
  we define the function $\rho $ by setting $\rho (  \tilde x, M): = \log \vert \tilde x M\vert$ for $( \tilde x, M) \in   \tilde {\mathbb X}\times \mathcal S$. This function satisfies the following cocycle property:  
  \[
  \rho (\tilde x, MN) = \rho (\tilde x, M)+\rho ( \tilde x \cdot M, N)
  \]
   for any $M, N \in\mathcal  S $ and $  \tilde x \in \tilde {\mathbb X}$.
   
    We achieve this paragraph with the  definition of a  {``lattice'' distribution} $\mu$  (see hypothesis {\bf P6}) in the context of products of random matrices.  It  may be stated  as follows:   the measure $\mu$ is {\it lattice}  if there exist  $t > 0, \epsilon \in [0, 2\pi[$ and a function $\psi: \mathbb X \to \mathbb R$   such that
	\begin{equation}\label{arith}
		\forall g \in T_\mu, \forall x \in \mathbb X, \quad \exp \left\{ it\rho(g, x)-i \epsilon +i (\psi(g\cdot x)-\psi(x)) \right\}=1,
	\end{equation}
	where $T_\mu$ is the closed sub-semigroup generated by the support of $\mu$.

\subsection{Exponential change of measures}

 
Under hypothesis {\bf P1}, for any $\theta \in [0, 1]$, we introduce the operator $P_\theta$  defined   by: for any bounded Borel function $\varphi: \tilde {\mathbb X}\to \mathbb C$ and $  \tilde x \in \tilde{\mathbb X}$,
\[
P_\theta \varphi( \tilde x):= \mathbb E[\vert   \tilde x M\vert^\theta \varphi( \tilde x \cdot M)].
\]
For any $n \geq 1$, it holds
\[
P^n_\theta \varphi( \tilde x):= \mathbb E[\vert  \tilde  x M_{0, n-1} \vert^\theta \varphi( \tilde x \cdot M_{0, n-1})].
\]
The operators $P_\theta$ are positive, they act continuously on the space $C(\tilde  {\mathbb X})$ of $\mathbb C$-valued continuous functions on $\tilde  {\mathbb X}$  endowed with the norm of uniform convergence $\vert \cdot \vert_\infty$. Under condition  {\bf P3},  their spectral radius on $C(\tilde  {\mathbb X})$ equals $\lambda(\theta)$.

 We denote by  $\mathcal B_{\theta}$ the space of   $\theta$-H\" older continuous functions  $\varphi: (\tilde  {\mathbb X}, d) \to \mathbb C$ such that
$
 m_{\theta} (\varphi):= \sup_{ \stackrel{ \tilde x,  \tilde y \in \tilde  {\mathbb X}}{  \tilde x \neq  \tilde y}}  {\vert \varphi( \tilde x )-\varphi( \tilde y)\vert \over  d( \tilde x,   \tilde y)^{\theta}}<+\infty.
$
Endowed with the norm $\vert \cdot \vert_{\theta}:= \vert \cdot \vert_\infty + m_{\theta} (\cdot),$ the space 
$(\mathcal B_{\theta}, \vert \cdot \vert_{\theta})$ is a $\mathbb C$-Banach space. Furthermore, for any $\varphi \in \mathcal B_{\theta}$ and $n \geq 1$, 
\begin{align*}
m_\theta(P^n_\theta \varphi) &\leq
\mathbb E[ c(M_{0, n-1})^\theta \vert M_{0, n-1}\vert ^\theta] \ m_\theta(\varphi)+ 2^\theta\  \mathbb E[\vert M_{0, n-1} \vert^\theta] \ \vert \varphi \vert_\infty\\
&
\leq
r^{n }\mathbb E[   \vert M_{0, n-1}\vert ^\theta] \ m_\theta(\varphi)+ 2^\theta\  \mathbb E[\vert M_{0, n-1} \vert^\theta] \ \vert \varphi \vert_\infty
\end{align*}
with  $r=  (c_{B})^\theta$ where $c_{B}$ is introduced in Properties \ref{propertiesHennion}. By hypothesis {\bf P3}, it holds $c_{B} \in [0, 1[$,  hence $r<1$. In other words, the operator $P_\theta$ satisfies the {\it Doeblin-Fortet} inequality  on $\mathcal B_\theta$: for any 
$\varphi \in \mathcal B_{\theta}$ and $n \geq 1$, 
\[
\vert P^n_\theta \varphi\vert_\theta \leq  r^n  \mathbb E[  \vert M_{0, n-1}\vert ^\theta] \vert  \varphi\vert_\theta+ C\times  \mathbb E[   \vert M_{0, n-1}\vert ^\theta] \vert \varphi\vert_\infty
\]
 for some constants $C>0$ and   $r \in ]0, 1[$. By 
  \cite{H93},  the operator $P_\theta$ is  quasicompact on $\mathcal B_\theta$, with spectral radius $\lambda(\theta)$. By \cite{GL} and \cite{BDGM},  the  real positive number $\lambda(\theta)$ is the unique eigenvalue of $P_\theta$ with modulus $\lambda(\theta)$, it is simple and there exist a unique strictly positive function $v_\theta \in \mathcal B_\theta$ and a unique probability measure $\nu_\theta$ on $\tilde  {\mathbb X}$  such that  $\nu_\theta (v_\theta)=1$   satisfying
  \[
  \nu_\theta P_\theta= \lambda(\theta) \nu_\theta \quad {\rm and}\quad  P_\theta v_\theta= \lambda(\theta) v_\theta.  \]
For any $\theta \in [0, 1]$, we thus introduce the new transition probability kernel $\bar{P}_\theta$ on $\tilde  {\mathbb X}$ defined by: for any bounded Borel function $\psi: \tilde  {\mathbb X} \to \mathbb C$  and any $\tilde x \in \tilde  {\mathbb X}$,
  \[ 
  \bar P_\theta \psi(\tilde x):= {1\over \lambda(\theta) v_\theta(\tilde x)}\mathbb E[ \vert \tilde xM\vert^\theta v_\theta(\tilde x \cdot M) \psi(\tilde x \cdot M)].
  \]
By the above, this Markovian operator   is quasicompact on $\tilde  {\mathbb X}$ with $1$ as  a simple and isolated eigenvalue and without any other eigenvalue of modulus 1. The powers of $\bar P_\theta$  are given by: for any $n \geq 1$, 
\[
  \bar P^n_\theta \psi(\tilde x):= {1\over \lambda(\theta)^n v_\theta(\tilde x)}\mathbb E[\vert \tilde x M_{0, n-1}\vert^\theta v_\theta(\tilde x \cdot M_{0, n-1}) \psi(\tilde x\cdot M_{0, n-1})].
  \]
 Let $(X^\theta_n)_{n \geq 0}$ be  the Markov chain on $\tilde  {\mathbb X}$ with transition operator $\bar P_\theta$. For $a \in \mathbb R^+$ fixed, any $\tilde x \in \tilde{\mathbb X}$ and any $n \geq 1$, we set
\[
S_0=a, \quad   S_n=S_n(\tilde x, a)= a+\ln \vert \tilde x M_{0, n-1}\vert.
\]   We associate to $\bar P_\theta$ the Markov  walk $(X^\theta_n, S_n)_{n \geq 0}$ on $\tilde  {\mathbb X}\times \mathbb R$ with transition operator $\tilde P_\theta$ defined by: for any bounded Borel function $\Psi: \tilde  {\mathbb X} \times \mathbb R\to \mathbb C$  and any $(\tilde x, a) \in \tilde  {\mathbb X}\times \mathbb R$,
  \[ 
  \tilde P_\theta \Psi(\tilde x, a):= {1\over \lambda(\theta) v_\theta(\tilde x)}\mathbb E[\vert \tilde x M\vert^\theta v_\theta(\tilde x \cdot M) \Psi(\tilde x \cdot M, a+\ln \vert \tilde x M\vert)].
  \]
In order to study the behavior of the $(X^\theta_n, S_n)_{n \geq 0}$, we introduce in a natural way the family of ``Fourier operators'' $(\bar P_{\theta, t})_{t \in \mathbb R}$ defined by:  for any bounded Borel function $\psi: \tilde  {\mathbb X} \to \mathbb C$,  for any $\theta \in [0, 1] $ and $t \in \mathbb R$, 
\begin{align*}
  \bar P_{\theta, t} \psi(\tilde x)&:= {1\over \lambda(\theta) v_\theta( \tilde x)}\mathbb E[\vert \tilde x M\vert^{\theta+it} v_\theta(\tilde x \cdot M) \psi(\tilde x \cdot M)].
\\
& =  {1\over \lambda(\theta) v_\theta( \tilde x)}\mathbb E[e^{(\theta+it)\ln\vert \tilde x M\vert} v_\theta(\tilde x \cdot M) \psi(\tilde x \cdot M)]
\end{align*}
In the sequel, it  is thus natural to consider, for any $\tilde x$ in $\tilde  {\mathbb X}$ and $a \in \mathbb R$, the  probability measures 
$ \mathbb P_{\tilde x, a}^\theta$ (with associated expectation $\mathbb E_{\tilde x, a}^\theta$) on $(\Omega, \mathcal F)$  whose restriction to the $\sigma$-algebras $\sigma(M_0, \ldots, M_{k-1}), k \geq 1$ is given by: for any positive Borel function $\Phi:  \mathcal S ^{k}\to \mathbb C$,
\begin{equation}\label{esperancetheta}
\mathbb E_{\tilde x, a}^\theta( \Phi(M_0, \ldots, M_{k-1})):= {1\over \lambda^k(\theta) v_\theta( \tilde x)}\mathbb E[\Phi(M_0, \ldots, M_{k-1})
e^{ \theta ( a+ \ln\vert \tilde x  M_{0, k-1})\vert} v_\theta(\tilde  x\cdot M_{0, k-1})]\
\end{equation}
and with the same conditional probability  measure $\mathbb P_{M_0=m_0, \ldots, M_{k-1}=m_{k-1}}$ as $\mathbb P$, for any $m_0, \ldots$, $m_{k-1}$ in  $ \mathcal S$.


 Now, we fix $\theta= \theta_\star$. The matrices $ M_n$ are iid with respect to the measure $\mathbb P$ but this property is no longer relevant with respect to $\mathbb P^{\theta_\star}_{\tilde x, a}$. Neverthelesss, due to the quasicompactness of the operators $\bar P_{\theta_\star}$ established above, the fluctuations of the process $(\tilde x M_{0, n-1} )$ with respect to $\mathbb P^{\theta_\star}_{\tilde x, a}$ may be controlled as in \cite{Pham2018}; the  Markov chain   $(M_n, \tilde x\cdot M_{0, n-1})_{n \geq 1}$ driven by the Markov operator $ P^{\theta_\star}$  falls within the scope of the article \cite{GLL2018} and  leads to the following  statement.

We denote by 
$$  \tau_{\tilde x, a}:= \min \{n \geq 1\mid S_n(\tilde x, a) \leq 0\}
$$
  the first moment when the sequence $(S_n(\tilde x, a))_{n \geq 0}$ enters the set 
$]-\infty, 0]$. 
  Modifying in a natural way the arguments used in \cite{Pham2018} one can conclude that under the conditions {\bf P1}--{\bf P6} the function $h^{\theta_\star}: \tilde{ \mathbb X}\times [0, +\infty[ \to [0, +\infty[$ specified by the equality
\[
h^{\theta_\star}(\tilde x, a) = \lim_{n \to +\infty} \mathbb E^{\theta_\star}(S_n(\tilde x, a), \tau_{\tilde x, a}>n)
\]
satisfies the property $\mathbb E^{\theta_\star}   ( h^{\theta_\star}(\tilde x\cdot M_0, S_1(\tilde x, a)); \tau_{\tilde x, a}>1)= h^{\theta_\star}(\tilde x, a)$. This function $h^{\theta_\star}$ satisfies the following properties (see Theorems 2.2 and 2.3 in \cite{GLL2018}).
\begin{prop} \label{fluctuationns} There exist  constants $C>1$ and $A>0$ such that   for any $x \in \mathbb X$ and $a>0$, 
\[
\mathbb P^{\theta_\star}(\tau_{\tilde x, a} >n) \sim {h^{\theta_\star}(\tilde x, a)\over \sqrt{2\pi n}} \quad   as  \quad n \to +\infty
\]
with
\[
\mathbb P^{\theta_\star} (\tau_{\tilde x, a} >n) \leq C\ {h^{\theta_\star}(\tilde x, a)\over \sqrt{ n}} \quad \ and \ \quad C^{-1}\max (1, a-A)\leq h^{\theta_\star}(\tilde x,a)\leq C (1+a).
\]
\end{prop}
In the sequel, we also need  a control of the quantities $\mathbb P^{\theta_\star} (S_n(\tilde x, a)\in [b, b+\ell], \tau _{\tilde x, a}>n)$ for any   interval $[b, b+\ell]$ included in $\mathbb R^+$. The following statement, in the spirit of Stone's theorem for classical random walks, is stated in \cite{GLL2017} for Markov walks over a finite state space and   is meaningful when the parameter $b$ is of  order $\sqrt{n}$. As announced in \cite{GLL2017},  this statement is valid in fact in  a more general setting, namely  when the underlying Markov chain is irreducible and  its transition operator  satisfies a spectral gap property.    When $\mu$ is non lattice the  spectral radius  in ${\mathcal B}_{\theta_\star}$ of  the  ``Fourier operators'' $\bar P_{{\theta_\star}, t}, t \in \mathbb R^*$,   is strictly less than $ 1$; this  is crucial in the proof of the following statement.

\begin{prop}\label{localconditionedStone}
Under hypotheses {\bf P1}--{\bf P6}, for any fixed $(\tilde x, a) \in \tilde{\mathbb X}\times \mathbb R, \ell>0$ and  uniformly in $b \in \mathbb R^+$, 
\[
\lim_{n \to +\infty} 
 \left(
 n\mathbb P^{\theta_\star} (\tau_{\tilde x, a}>n, S_n({\tilde x, a})\in [b, b+\ell]) - {2\  \ell\ h^{\theta_\star}(\tilde x, a)\over \sqrt{2\pi \sigma^2} }\varphi_+\left({b\over \sigma \sqrt{n} }\right)
 \right)=0,
\]
where $\varphi_+(t)= t e^{-t^2/2}{\bf 1}_{\mathbb R^+}(t)$ is the Rayleigh density.
\end{prop}
This is now  of interest to introduce the following notations.    For any $y \in \mathbb X$ and $b >0$, we set 
$$
\tilde S_0( y, b)=b \quad \text{and} \quad \tilde S_n( y, b):= b-\ln \vert M_{n-1, 0}y\vert  \quad \text{for} \quad    n \geq 1.
$$
We
denote by 
$$ \tilde \tau_{y, b}:= \min \{n \geq 1\mid \tilde S_n(y, b) \leq 0\}
$$
  the first moment when the sequence $(\tilde S_n(y, b))_{n \geq 0}$ enters the set 
$]-\infty, 0]$.  Propositions  \ref{fluctuationns}  and  \ref{localconditionedStone}  also hold  for the process $(M_{n-1, 0}\cdot y, \tilde{S}_n(y, b))_{n \geq 1}$  and the stopping time $\tilde \tau_{y, b}$, we denote by $\tilde h^{\theta_\star}$ the    function analogous to $h^{\theta_\star}$ which appears in the corresponding  statements for the process $(M_{n-1, 0}\cdot y, \tilde{S}_n(y, b))_{n \geq 1}$. 

 In the sequel, we set $\Delta:= \ln \delta$ where $\delta$ is the constant which appears in Lemma \ref{keylemma}; the proof  is detailed in \cite{PP2023}.

 Propositions  \ref{fluctuationns}  and    \ref{localconditionedStone}   yield  to   following  statement. 
\begin{coro}\label{theolocalPthetastar}
		 Assume  hypotheses {\bf P1}--{\bf P6} hold.  Then, there exists a   positive constant    $  C>0$   such that  for any $ \tilde x \in \tilde {\mathbb X}, a, b\geq 0$  and  any   $\ell>0,$  
		\begin{equation}\label{3/2upper}
			 \mathbb P ^{\theta_\star} (\tau_{\tilde x, a}  >n, S_n(\tilde x, a)  \in  [b, b+\ell])\leq {C\over n^{3/2}}\  h^{\theta_\star}(\tilde x, a)\ {\tilde h}^{\theta_\star}(\tilde x, b) \ \ell.
		\end{equation}
 Furthermore,   there exist  constants $c>0$ and $\ell_0$ such that for $\ell>  \ell_0$  and $b \geq {B}$,  
			\begin{equation}\label{3/2lower}
				\liminf_{n \to +\infty}n^{3/2} \mathbb P_{\tilde x, a}^{\theta_\star}( \tau_{\tilde x, a} >n, S_n(\tilde x, a)  \in  [b, b+\ell])\geq c\  h^{\theta_\star}(\tilde x, a)\ {\tilde h}^{\theta_\star}(\tilde x, b)\  \ell.
			\end{equation}
	\end{coro}
Proof of corollary \ref{theolocalPthetastar}. 
When studying fluctuations of  random walks  $(S_n)_{n \geq 1} $ with iid increments $Y_k$ on $\mathbb R^p, p\geq 1$, we often reverse  time as follows.  For any $1\leq k \leq n$, the   random variables   $S_n-S_k= Y_{k+1}+\ldots +Y_n$  and $S_{n-k}=Y_1+\ldots +Y_{n-k}$ have the same distribution. In the case of products of random matrices, even when the $M_i$ are iid (which is not the case under  $\mathbb P^{\theta_\star}$), the   cocycle property $S_n(\tilde x)=S_n(\tilde x, 0) = \ln \vert \tilde x M_{0, n-1} \vert = S_k(\tilde x)+S_{n-k}(\tilde X_k),$ with $\tilde X_k= \tilde x \cdot M_{0, k-1}$,    decomposes  the sum $S_n(\tilde x)$ into two terms which are not independent. Hence, the  same argument  as in the iid  case cannot be applied directly. The fact that the  matrices $M_0, M_1, \ldots $ belong to $\mathcal S_{B}$ helps here:  for any  $x, y \in \mathbb X$,  we can compare the distribution of $S_n(\tilde x)-S_k(\tilde x)$  to   the one of $\ln \vert   M_{ n-k-1, 0}y\vert  $   (notice here that   the  non commutativity of the product of matrices forces us to  consider in this last quantity  the right  linear  action     of the matrices $M_{0, n-k-1}$). It is the strategy that we apply to obtain the following result (see lemma 2.3  in \cite{PP2023}) which is central in the sequel.

\begin{lemma}  \label{reverse} 
		For any $  x, y \in \mathbb X, a, b \geq   0$ and $ \ell>0$, 
		\begin{equation*}\label{reversingineq1}
			\mathbb P (\tau_{\tilde x, a} >n,   S_n(\tilde x, a) \in [b, b+\ell])\leq \mathbb P(\tilde \tau_{  y, \ b+\ell+{B}} >n,   \tilde S_n(y, \ b+\ell+{B})  \in [a, a+\ell+2{B}]).
		\end{equation*}
		Similarly,  for $a\geq  \ell>2{B} >0$  and $b \geq \Delta$,
		\begin{equation*}\label{reversingineq2}
			\mathbb P(\tau_{\tilde x, a} >n,   S_n(\tilde x, a)\in [b, b+\ell])\geq \mathbb P(\tilde \tau_{ y, \ b-{B}} >n,   \tilde S_n(y, \ b-B) \in [a-\ell, a-2{B}]).
		\end{equation*}
\end{lemma}
This allows us to apply the same strategy as in \cite{PP2023} to prove Corollary \ref{theolocalPthetastar}. Let us propose a sketch of the argument.

\noindent {Proof.}
Inequality (\ref{3/2upper}) is proved in \cite{LPP2021} Corollary 3.7.  The proof of  the lower bound  $(\ref{3/2lower})$  is based on Proposition \ref{localconditionedStone}. We set $m=\lfloor n/2 \rfloor$ ; by
  the Markov property and  Lemma \ref{reverse},
\begin{align}\label{zeudg}
	&\mathbb P^{\theta_\star}  ( \tau_{\tilde x, a} >n,\   S_n(\tilde x, a)\in  [ b, b+\ell ])\notag \\
	&\geq
	\mathbb P^{\theta_\star}( \tau_{\tilde x, a} >n,  S_{m}(\tilde x, a)\in [\sqrt n, \sqrt{2n}],  S_n(\tilde x, a)\in  [ b, b+\ell ])
	\notag \\
	&\geq 
	\sum_{\stackrel{k\in \mathbb N}{\sqrt{n}\leq \ k\  \leq \sqrt{2n}-1}} 
	\mathbb P^{\theta_\star}\Bigl( \tau_{\tilde x, a} >n,  k \leq S_{m}(\tilde x, a) \leq  k+1, \ b \leq  S_n(\tilde x, a)\leq b+\ell \Bigr)
	\notag \\
	&\geq 
	\sum_{\stackrel{k\in \mathbb N}{\sqrt{n}\leq \ k \ \leq \sqrt{2n}- 1}} 
	\int_{\tilde {\mathbb X}\times [k,  k+1 ]} \mathbb P^{\theta_\star} ( \tau_{\tilde x, a} >m, (\tilde x\cdot M_{0, m-1}, S_m(\tilde x, a)) \in {\rm d}\tilde x'{\rm d}a')\notag \\
	&  \qquad \qquad \qquad   \mathbb P^{\theta_\star}(\tau_{{\tilde x}', a'} >n-m, b  \leq  S_{n-m}({\tilde x}', a') \leq b +\ell )
	\notag \\
	&
	\geq 
	\sum_{\stackrel{k\in \mathbb N}{\sqrt{n}\leq \ k \ \leq \sqrt{2n}- 1}} 
	\int_{\tilde{\mathbb X}\times [k,  k+1 ]} \mathbb P^{\theta_\star}\Bigl( \tau_{\tilde x, a} >m, (\tilde x\cdot M_{0, m-1}, S_m(\tilde x, a))  \in {\rm d}\tilde x'{\rm d}a')\notag \\
	&  \qquad\qquad \qquad  \mathbb P^{\theta_\star}(\tilde \tau_{  x, b-\Delta} >n-m, a'-\ell   \leq \tilde S_{n-m}(  x, b-\Delta) \leq  a'-2 \Delta)
	\notag \\
	&
	\geq 
	\sum_{\stackrel{k\in \mathbb N}{\sqrt{n}\leq \ k \leq \sqrt{2n}- 1}} 
	\mathbb P^{\theta_\star}\Bigl( \tau_{\tilde x, a} >m, k  \leq S_m(\tilde x, a)\leq  k+1 \Bigr) \notag\\
	&  \qquad\qquad \qquad  \mathbb P^{\theta_\star}\left(\tilde \tau_{x, b-\Delta}>n-m, k +1 -\ell \leq \tilde S_{n-m}(x, b-\Delta) \leq  k -2 \Delta\right)
\end{align}  
By Proposition  \ref{localconditionedStone},  there exists a    constant  $C_0>0$ such that when   $\sqrt{n}\leq\  k \displaystyle \ \leq \sqrt{2n}- 1$,
\[
\liminf_{n \to +\infty}n \mathbb P^{\theta_\star}\Bigl( \tau_{\tilde x, a} >m, k  \leq S_m(\tilde x, a)\leq  k+1 \Bigr)\geq C_0 \  h^{\theta_\star} (\tilde x, a)
\]
and
\[
\liminf_{n \to +\infty}\mathbb P^{\theta_\star}\left(\tilde \tau_{x, b-\Delta}>n-m, k+1-\ell  \leq \tilde S_{n-m}(x, b-\Delta) \leq  k -2 \Delta\right)\geq C_0 (\ell -2 \Delta-1) \   \underbrace{{\tilde h}^{\theta_\star} (\tilde x, b-\Delta)}_{\succeq  \ {\tilde h}^{\theta_\star} (\tilde x, b)  \ }.
\]
Hence, inequality (\ref{zeudg}) yields, for $n$ large enough, 
\begin{align*}
	n^2 \mathbb P^{\theta_\star}( \tau_{\tilde x, a} >n,  S_n(\tilde x, a)\in  [ b, b+\ell ])\geq {C_0^2\over 2}   (\sqrt{2n}-\sqrt{n})  (\ell - 2 \Delta-1)\ h^{\theta_\star} (\tilde x, a)\  {\tilde h}^{\theta_\star} (\tilde x, b),
\end{align*}
which   implies, for such $n$,  
\[
\mathbb P^{\theta_\star}( \tau >n,  S_n\in  [ b, b+\ell ])
\succeq   {h^{\theta_\star} (\tilde x, a) \  {\tilde h}^{\theta_\star} (\tilde x, b)  \over n^{3/2}}(\ell - 2 \Delta-1).  \]
This achieves the proof of the lower bound (\ref{3/2lower}) taking $\ell_0:= 4\Delta+2$.

\rightline{$\Box$}

As a direct consequence,  we obtain the following  ``rough local limit theorem'' for the process $(S_n(\tilde x, a))_{n \geq0}$  and its minimum under the probability $\mathbb P $.
\begin{coro}\label{corolocalnoncentered}
		Under hypotheses {\bf P1}--{\bf P6}, there exist  positive constants  $c, C$ and $\ell_0$ such that  for any $a \geq \ell \geq \ell_0,  b \geq \Delta$     and any $  x  \in  \mathbb X$,   
		\begin{equation}\label{localnoncenteredupper}
		 \mathbb P (\tau_{\tilde x, a}  >n,   S_n(\tilde x, a) \in  [b, b+\ell])  
 \leq  C  {\rho_\star^n\over n^{3/2}}  \ h^{\theta_\star} (\tilde x, a)\  {\tilde h}^{\theta_\star} (\tilde x, b)\ e^{\theta_\star(a-b)} \ \ \ell 
	\end{equation}
and 		
\begin{equation}\label{localnoncenteredlower}
				\liminf_{n \to +\infty} {n^{3/2}\over \rho_\star^n}   \mathbb P (\tau_{\tilde x, a}  >n,   S_n(\tilde x, a) \in  [b, b+\ell])   \geq c\  h^{\theta_\star} (\tilde x, a)\  {\tilde h}^{\theta_\star} (\tilde x, b)\ e^{\theta_\star(a-b-\ell)}\ell  .
			\end{equation}
%
%
Similarly,  
\begin{equation}\label{minimumnoncentered}
\mathbb P  (\tau_{\tilde x, a}>n)  
		\leq  C  {\rho_\star^n\over n^{3/2}} \   e^{\theta_\star a} \ h^{\theta_\star} (\tilde x, a)\quad {\it and} \quad 
\liminf_{n \to +\infty}  {n^{3/2}\over \rho_\star^n}  \mathbb P  (\tau_{\tilde x, a}>n)    \geq   c\  e^{\theta_\star a} \  h^{\theta_\star} (\tilde x, a). 
\end{equation}
\end{coro}
Proof. By using definition of $\mathbb E^{\theta_\star}$ (see 
(\ref{esperancetheta}))  and the fact that the function $v_{\theta_\star}$ is non negative and continuous on $\tilde{\mathbb X}$, we may write 
\begin{align*} 
		 \mathbb P (\tau_{\tilde x, a}  >n,  S_n(\tilde x, a) \in  [b, b+\ell]) \asymp
		 \rho_\star^n\ e^{\theta_\star a} \ \mathbb E^{\theta_\star}  [\tau_{\tilde x, a}  >n,  e^{-\theta_\star S_n(\tilde x, a)} {\bf 1}_{ [b, b+\ell]}(S_n(\tilde x, a))] 
		  		\end{align*}
				with $e^{-\theta_\star(b+\ell)}\ {\bf 1}_{ [b, b+\ell]}(S_n(\tilde x, a))\leq  e^{-\theta_\star S_n(\tilde x, a)} \ {\bf 1}_{ [b, b+\ell]}(S_n(\tilde x, a))\leq e^{-\theta_\star b}\ {\bf 1}_{ [b, b+\ell]}(S_n(\tilde x, a))$.
				Similarly $\displaystyle \mathbb P (\tau_{\tilde x, a}  >n) \asymp
		 \rho_\star^n\ e^{\theta_\star a} \ \mathbb E^{\theta_\star} [\tau_{\tilde x, a}  >n,  e^{-\theta_\star S_n(\tilde x, a)} ] $ with 
		 \[
		 \mathbb E^{\theta_\star} [\tau_{\tilde x, a}  >n,  e^{-\theta_\star S_n(\tilde x, a)} ]\asymp \sum_{k\geq 0}
		 e^{-\theta_\star k} \mathbb P^{\theta_\star}(\tau_{\tilde x, a}  >n, S_n(\tilde x, a)\in [k, k+1[).\]
				
				\noindent Inequalities (\ref{localnoncenteredupper}),  (\ref{localnoncenteredlower})   and (\ref{minimumnoncentered}) follow, by applying Corollary \ref{theolocalPthetastar}.
				\vspace{2mm}
				
\rightline{$\Box$}

%
%

	\section{Proof of the main theorem}
		 	 \subsection{ Expression of non extinction  probability}  
			   For  every realization of the environmental sequence ${\bf f}= (f_0, f_1, \ldots)$ and $0\leq k \leq n-1$, we define
 $$F_{k, n-1}= f_k\circ \ldots \circ f_{n-1}\quad {\rm and} \quad  F_{ n, n-1}= {\rm Id}.
 $$
 For $1\leq i \leq p$, we set  $F^{(i)}_{k, n-1}= f_k^{(i)}\circ \ldots \circ f_{n-1}$.
From the definition of the process $(Z_n)_{n \geq 0}$, it holds
 \[
 \mathbb E[s^{Z_n}\mid Z_0= \tilde e_i,  f_0, \ldots f_{n-1}]= f_0^{(i)}\circ \ldots \circ f_{n-1}(s)= F_{0, n-1}^{(i)}(s) 
 \]
 so that
 \[
\mathbb E [q_n^{(i)}]=  \mathbb P(\vert Z_n\vert >0\mid Z_0= \tilde e_i)= \mathbb E[1-F_{0, n-1}^{(i)}({\bf 0})]. \]
 For a generating function  $f$ with corresponding mean matrix $M$ and a matrix ${\bf a}= ({\bf a}(k, \ell))_{1\leq k, \ell \leq p}$ with positive entries, we set, for   $s \in [0, 1]^p$, 
 $$
 \psi_{f, {\bf a}}(s):= {\vert {\bf a}\vert \over  \vert {\bf a}({\bf 1}-f(s))\vert}- {\vert {\bf a}\vert \over  \vert  {\bf a}M({\bf 1}-s)\vert}.
 $$
We fix $i \in \{1, \ldots, p\}$.
 Let ${\bf a}^{(i)}   $  be the matrix with ${\bf a}^{(i)}(i, i)=1$ and $ {\bf a}^{(i)}(k, \ell)=0$ for all $(k, \ell)\neq (i, i)$.  Using the definition of the functions   $\psi_{f, {\bf a}}$, we write for $n \geq 1$,
\begin{align*}
{1\over  1-F_{0, n-1}^{(i)}(s)}
 =
{1\over  \vert  {\bf a}^{(i)} M_{0, n-1}({\bf 1}-s)\vert} + \sum_{k=0}^{n-2} {1\over \vert {\bf a}^{(i)}M_{0, k}\vert}
\psi_{f_k, {\bf a}^{(i)}M_{0, k-1}} (F_{k+1, n-1}(s)),
\end{align*}
where, as previously $M_{0, k}= M_0\ldots M_k$ for any $k \geq 0$. Consequently,
\begin{align*}  
\mathbb  E[q_n^{(i)}] &= \mathbb E[1-F^{(i)}_{0, n-1}({\bf 0})]\ \\
&= \mathbb E\left[\left(
{1\over  \vert  {\bf a}^{(i)} M_{0, n-1}  \vert} + \sum_{k=0}^{n-2} {1\over \vert {\bf a}^{(i)}M_{0, k}\vert}\underbrace{\psi_{f_k, {\bf a}^{(i)}M_{0,k-1}} (F_{k+1, n-1}({\bf 0}))}_{\eta_{k, n-1} } 
\right)^{-1}\right]
 \\
&= \mathbb E\left[\left(
   \sum_{k=0}^{n-1} {1\over \vert {\bf a}^{(i)}M_{0, k}\vert}\underbrace{\psi_{f_k, {\bf a}^{(i)}M_{0,k-1}} (F_{k+1, n-1}({\bf 0}))}_{\eta_{k, n-1} } 
\right)^{-1}\right] 
\end{align*}
where the random variable $\eta_{k, n-1}$ is defined as above with the convention $\eta_{0, n-1}=1$. By Lemma  3 in \cite{DV2018}, it holds $0\leq \eta_{k, n-1} \leq {B} p^2 \sum_{i=1}^p \vert B_k^{(i)}\vert / \vert M_k\vert^2$. 

In the sequel, we also need  a control of the lower bound of the  random variables $\eta_{k, n-1}$. This is why we introduce the restrictive assumption {\bf P7}, which readily implies,  by a straightforward computation, that the random variables are bounded from below and above by non negative constants (see Proposition 2.1 and Lemma 3.1 in \cite{DHKP}). Therefore, under this additive assumption, it holds 
\begin{align} \label{extinctionprob}
\mathbb  E[q_n^{(i)}] \asymp  \mathbb E\left[ {1\over 
\vert  M_{0, 0}\vert ^{-1}+\ldots + \vert M_{0, n-1}\vert^{-1}} \right].
\end{align}
 This is not an exact formula but it is sufficient in our context since we only have rough estimates in Corollary  \ref{corolocalnoncentered}. We cannot improve our  result as long as we do not have a precise estimate in the local theorem for the norms of products of random matrices conditioned to remain strictly greater than 1 until time  $n$.

 \subsection{The change of measure $ \widehat {\mathbb P}_{\tilde x, a}$}
%
In the case when the Lyapunov exponent of the sequence $(M_n)_{n \geq 0}$ equals $0$, it is natural to introduce the following 
new probability measure $ \widehat {\mathbb P}_{\tilde x, a}$  on the canonical path space $ ((\mathbb{X} \times \mathbb{R})^{\otimes \mathbb{N}}, \sigma (X_n, S_n: n \geq 0), \theta)$   of  the Markov chain $(X_n, S_n)_{n \geq 0}$ $^($\footnote{$\theta$ denotes   the shift operator on $(\mathbb{X} \times \mathbb{R})^{\otimes \mathbb{N}}$  defined by $\theta \Bigl((x_k, s_k)_{k\geq 0}\Bigr)=  (x_{k+1}, s_{k+1})_{k\geq 0} $ for any $(x_k, s_k)_{k\geq 0}$  in $(\mathbb{X} \times \mathbb{R})^{\otimes \mathbb{N}}$}$^)$, whose associated expectation $ \widehat {\mathbb E}_{\tilde x, a}$ is   characterized by the following property: for any bounded and non negative Borel function $\varphi$ on $(\tilde {\mathbb{X}}\times \mathbb{R})^{k+1}$ with compact support,
\begin{align*}   \label{lcaehjrzgbc}
\widehat {\mathbb{E}}_{\tilde x,a}[\varphi (X_0, S_0,\ldots,X_k, S_k)] &:= 
 \lim_{n \to {+\infty} } \mathbb{E}_{\tilde x,a} [\varphi (X_0, \ldots,  S_k)\mid \tau >n], 
\end{align*}
when the limit exists.
The fact that the Lyapunov   exponent of the sequence $(M_n)_{n \geq 0}$ equals $0$ implies that the probability $ \mathbb P_{\tilde x, a}(\tau>n)$    decreases  towards $0$ as $1/\sqrt{n}$ (see \cite{Pham2018}) and the expectation $\widehat {\mathbb E}_{\tilde x,a}$ is given explicitely as follows
\[
 \widehat {\mathbb{E}}_{\tilde x,a}[\varphi (X_0, S_0,\ldots,X_k, S_k)]={1\over h(\tilde x, a)}{\mathbb{E} _{\tilde x,a}}[\varphi ({X_0,\ldots, S_k})h(X_k, {S_k}),\, \tau>k],
 \]
 for some non negative  function  $h$ (see \cite{LPP2018} for the details). Furthermore,   if $(Y_k)_{k \geq 0}$  is a sequence of  bounded real-valued random variables,  adapted to the filtration  $(\mathcal F_k)_{k \geq 0}$ and which  converges  in $\mathbb L^1(\widehat {\mathbb P} _{\tilde x, a})$ to some  random variable   $Y_\infty$, then 
		\begin{equation} \label{oaezrtufnqj}
		\lim_{n \to +\infty } \mathbb E_{\tilde x, a} \big[Y_n \ \vert \ \tau >n\big] = \widehat {\mathbb E}_{\tilde x, a}[Y_\infty].
		\end{equation}
This last property yields to  the speed of convergence to $0$ of the extinction probability 
in a rather luminous way \cite{LPP2018},\cite{LPP2021}.

 In the present case,  a similar construction does exist under the probability $\mathbb P^{\theta_\star}$;  in this case, the sequence $(M_n)_{n \geq0}$ is not iid but the choice of the parameter $\theta_\star$ is done in such a way the corresponding Lyapunov exponent $\Lambda'(\theta_\star) $ equals $0$.  Therefore, we introduce the following family of probability measures $\widehat{\mathbb P}^{\theta_\star}_{\tilde x, a}$  (with corresponding expectation $\widehat{\mathbb E}^{\theta_\star}_{\tilde x, a}$) whose restriction to the $\sigma$-algebras
 $\sigma(X_0, S_0, \ldots, X_k, S_k), k\geq 0$ are defined as follows: for any bounded and non negative Borel function $\varphi$ on $(\mathbb{X}\times \mathbb{R})^{k+1}$ with compact support
\begin{align}   \label{se;jrtbv}
\widehat {\mathbb{E}}_{\tilde x,a}^{\theta_\star}[\varphi &(X_0, S_0,\ldots,X_k, S_k)] \notag\\
& := 
{1\over h^{\theta_\star}(\tilde x, a)}{\mathbb{E}^{\theta_\star} _{\tilde x,a}}[\varphi ({X_0,\ldots, S_k})h^{\theta_\star}(X_k, {S_k}),\, \tau>k]\notag
\\
&= 
{\mathbb E _{\tilde x, a}\left[\varphi (X_0,\ldots, S_k)h^{\theta_\star}(X_k, {S_k})e^{\theta_\star S_k}
 v_{\theta_\star}(X_k),\, \tau>k\right]\over \rho_\star^k \ h^{\theta_\star}(\tilde x, a) \ v_{\theta_\star}(\tilde  x)}
\end{align}
The rough estimate of the quantity $ \mathbb P_{\tilde x, a}(\tau >n)$     given by (\ref{minimumnoncentered})  immediately yields   
\begin{align} \label{hsqdkhfkla}
\widehat {\mathbb{E}}_{\tilde x,a}^{\theta_\star}[\varphi (X_0, S_0,\ldots,X_k, S_k)] 
 &\asymp 
 {\mathbb{E}_{\tilde x, a}\left[\varphi (X_0,\ldots, S_k)h^{\theta_\star}(X_k, {S_k})e^{\theta_\star S_k}
,\, \tau>k\right]\over \rho_\star^k \ h^{\theta_\star}(\tilde x, a)}. 
\end{align}
Let us conclude this paragraph noticing that  the property (\ref{oaezrtufnqj}) does not hold anymore, even in a weaker form,   as the reader can see by following the demonstration of lemma 4.1 in  \cite{LPP2021}.
The strategy to get the speed of convergence towards $0$  of the survival probability is thus different and uses a clever decomposition of  the denominator of the right hand side in (\ref{extinctionprob}), inspired by \cite{LP1997} and \cite{GKV2002}.
 


\subsection{Proof of Theorem \ref{theoprincipal}}
 
We set here $m_k:= \min (\vert M_{0, 0}\vert, \ldots, \vert M_{0, k-1}\vert)= \min (\vert \tilde {\bf 1}M_{0, 0}\vert, \ldots, \vert \tilde {\bf 1}M_{0, k-1}\vert)$ for any $k \geq 1$.

\noindent Hence $(m_n\geq e^{-a})= (\tau_{\tilde {\bf 1}, a}>n)$ for any $a\geq 0$ and (\ref{localnoncenteredlower}) and (\ref{minimumnoncentered}) may be restated in particular as follows:  for $\ell > \ell_0$,  
\begin{align}\label{localminimumliminf}
		\liminf_{n \to +\infty} {n^{3/2}\over\rho_\star^n  }  \ \mathbb P (\tau_{\tilde x, a}  >n,   S_n(\tilde x, a) \in  [b, b+\ell]) 
\geq c  h^{\theta_\star} (\tilde x, a)\  {\tilde h}^{\theta_\star} (\tilde x, b) \ e^{\theta_\star(a-b-\ell)}\ \ell
		\end{align}
and 
\begin{align}
\label{minimumbis}
		  \mathbb P  (m_n \geq e^{-a})
			\leq \  C \ {\rho_\star^n\over n^{3/2}} \   e^{\theta_\star a}  h^{\theta_\star} (\tilde {\bf 1}, a),
\end{align}
where $c$ and $C$ are two strictly positive constants.
\subsubsection{Proof of the upper bound} 
It holds 
\begin{align*}
\mathbb P (\vert Z_n \vert >0\mid Z_0= \tilde{e_i})&= \mathbb E [q_n^{(i)}]
\\
&\leq   \mathbb P (m_n \geq 1)+ \mathbb E [q_n^{(i)}; m_n <1]
\\
&\preceq  \mathbb P (m_n \geq 1)+ \mathbb E[m_n; m_n <1]
\\
&\leq  \mathbb P (m_n \geq 1)+  \sum_{j =1}^ {+\infty} {{e^{ - j+1}}\mathbb{P}(e^{-j} \leq   m_n < e^{ - j+1})}  \\
&\leq  \mathbb P (m_n \geq 1)+  \sum_{j = 1}^ {+\infty}  e^{ - j+1}  \mathbb P( m_n\geq  e^{ - j}) 
\\
&\preceq {\rho_\star ^n\over n^{3/2}}+   {\rho_\star^n\over n^{3/2}}\sum_{j=1}^ {+\infty}  h^{\theta_\star} (\tilde {\bf 1}, j) \ e^{ (\theta_\star-1) j}
\quad {\rm by} \ (\ref{minimumbis}) ,
\\
&  \preceq \quad {\rho_\star^n\over n^{3/2}} \quad  \text{since} \quad  h^{\theta_\star} (\tilde {\bf 1}, j)\preceq 1+j\quad  \text{and} \quad \theta_\star <1.
\end{align*}
The upper bound is established.

\subsubsection{Proof of the lower  bound}

We fix $0< k < n/2$ and decompose 
$\displaystyle \sum_{\ell=0}^{n-1} {1\over \vert  M_{0, \ell}\ \vert}     $
 as ${\bf A}+{\bf B}+{\bf C}$ where
\[
{\bf A}  :=  \sum_{\ell=0}^{k-1} \vert M_{0, \ell}\vert ^{-1}, \quad  {\bf B} :=  \sum_{\ell=k}^{n-k-1} \vert  M_{0, \ell}\vert ^{-1}  
\quad{\rm and} \quad   {\bf C} :=    \sum_{\ell=n-k}^{n-1} \vert  M_{0, \ell}\vert ^{-1}   .
\]
It holds 
\begin{align*}
\mathbb P (\vert Z_n \vert >0\mid Z_0= \tilde{e_i})&= \mathbb E [q_n^{(i)}]\\
&\geq   \mathbb E \left[{1\over 1+{\bf A}+{\bf B}+{\bf C}}; m_n \geq 1\right]\\
&\geq  \mathbb E \left[{1\over1+ {\bf A}+{\bf C}}; m_n \geq 1\right]- \mathbb E \left[{\bf B}; m_n \geq 1\right].
 \end{align*}
  
\noindent  {\bf First step: control of the term} $\displaystyle   \mathbb E \left[{1\over1+ {\bf A}+{\bf C}}; m_n \geq 1\right]$ 
 
 \noindent It holds
  \begin{align*}
  \mathbb E \left[{1\over {\bf A}+{\bf C}}; m_n \geq 1\right] &
  \geq    \mathbb E [{1\over 1+{\bf A} + \sum_{\ell=n-k}^{n-1} \vert M_{0, \ell}\vert ^{-1} } ; 
   \\
   &   \qquad    m_{n-k} \geq 1, \vert  M_{n-k, n-1}\vert \leq 1, \ldots, \vert  M_{n-1, n-1}\vert \leq 1  ]\\
   & \quad {\rm since} \ [m_{n-k}\geq 1,  \vert  M_{n-k, n-1}\vert \leq 1, \ldots, \vert  M_{n-1, n-1}\vert \leq 1]\subset [m_n \geq 1]
   \\
   &\geq    \mathbb E [{1\over 2+{\bf A} + \sum_{\ell=n-k}^{n-1} \vert M_{\ell, n-1}\vert } ; 
   \\
   &   \qquad    m_{n-k} \geq 1, \vert  M_{n-k, n-1}\vert \leq 1, \ldots, \vert  M_{n-1, n-1}\vert \leq 1  ]
   \\
   & \quad {\rm since} \  \vert M_{\ell+1, n-1}\vert \geq \vert M_{0, \ell} \vert^{-1}, n-k\leq \ell \leq n-2 \ {\rm and } \ \vert M_{0, n-1}\vert^{-1} \leq 1 \\
   & \quad  {\rm on \ the \ event } \ [\vert M_{0, n-1}\vert \geq 1]   \\
  &\geq
   \mathbb E \left[{1\over  2+{\bf A} + \sum_{\ell=n-k}^{n-1} \vert  M_{\ell, n-1}\vert  }; \right.
 \\
   & \qquad     \left. m_{n-k} \geq 1, 1\leq \vert M_{0, n-k -1}\vert \leq K,  {1\over K} \leq \vert M_{n-k, n-1}\vert   \leq 1,  \ldots, \vert  M_{n-1, n-1}\vert \leq 1\right]
   \\
   & \quad {\rm for \ any \ constant \ } K>1,
   \\ 
      &\geq
   \mathbb E \left[{1\over  2+{\bf A} + \sum_{\ell=1}^{k} \vert   M'_{\ell, 1}\vert  }; \right.
   \\
   &   \qquad    \left. m_{n-k} \geq 1, 1\leq \vert M_{0, n-k -1}\vert \leq K,  {1\over K} \leq \vert  M'_{k, 1}\vert   \leq 1, \ldots,  \vert  M'_{ 1,  1}\vert\leq 1\right]\\ 
   & \quad {\rm by \ setting} \  (M_{n-k}, \ldots, M_{n-1})= (M'_{k}, \ldots, M'_1) \  {\rm and} \quad 
  M'_{\ell, 1}= M'_\ell \ldots M'_1\\
 &=
   \mathbb E  \left[{1\over  2+{\bf A}+{\bf C}'}; m_{n-k} \geq 1, 1\leq \vert M_{0, n-k -1}\vert \leq K,  {1\over K} \leq \vert  M'_{k, 1}\vert   \leq 1,  {\bf M}'_{k}\leq 1\right]\\ 
   &\quad  {\rm with } \ {\bf C}':= \sum_{\ell=1}^{k} \vert   M'_{\ell, 1}\vert   \  {\rm and} \ 
 {\bf M}'_{k}:= \max (\vert  M'_{1, 1}\vert, \ldots, \vert  M'_{ k, 1}\vert)
\\
  &
  \geq {1\over 2} \mathbb E  \left[{1\over  1+{\bf A}}\times {1\over 1+{\bf C}'}; m_{n-k} \geq 1, 1\leq \vert M_{0, n-k -1}\vert \leq K, \right.
  \\
  & \qquad \qquad \qquad \qquad \qquad \qquad \qquad \qquad  \left. {1\over K} \leq \vert  M'_{k, 1}\vert   \leq 1,  {\bf M}'_{k}\leq 1 \right]
  \end{align*}
\begin{align*}
 &
 \geq  {1\over 2}\mathbb E \left[{1\over  1+{\bf A} }; m_{n-k} \geq 1, 1\leq \vert M_{0, n-k -1}\vert \leq K\right] \\
 &  \qquad \qquad \qquad \qquad \qquad \qquad \qquad \qquad  \times \mathbb E \left[{1\over  1+{\bf C} '};   {1\over K} \leq \vert  M'_{k, 1}\vert   \leq 1, {\bf M}'_{k} \leq 1\right]
 \\
 &  \quad {\rm by \ independence\ of\ the\ random \ variables \  } \ {1\over 1+{\bf A} } {\bf 1}_{[m_{n-k} \geq 1, 1\leq \vert M_{0, n-k -1}\vert \leq K]} {\rm \ and}\ \\
 & 
 \quad  {1\over 1+{\bf C} '} {\bf 1}_{[{1\over K} \vert  M'_{k, 1}\vert   \leq 1, {\bf M}'_{k} \leq  1]}.
  \end{align*}
Now, we assume  $\ln K\geq \ell_0$, where $\ell_0$ is defined in Lemma \ref{reverse}; it holds, on the one hand, 
\begin{align*}
\liminf_{n \to +\infty} {n^{3/2} \over \rho_\star^{n-k}}\mathbb E &\left[{1\over  1+{\bf A} };   m_{n-k} \geq 1, 1\leq \vert M_{0, n-k -1}\vert \leq K\right] 
\\
& \geq \liminf_{n \to +\infty} {n^{3/2} \over \rho_\star^{n-k}}
\mathbb E \left[{1\over  1+{\bf A} };  m_{k} \geq 1, \mathbb P_{\tilde x_k, a_k}  (
m_{n-2k}\circ \theta^k\geq 1, 1\leq \vert M_{0, n-2k  }\circ \theta^k\vert \leq K)\right]
\\
& {\rm with\ }  \tilde x_k=\tilde{\bf 1}\cdot M_{0, k-1} \ {\rm and} \ a_k= \ln\vert M_{0, k-1}\vert
\\
& \geq \rho_\star^{-k}
\mathbb E \left[{1\over  1+{\bf A} };  m_{k} \geq 1,\right.
\\
& \qquad \qquad \qquad \underbrace{ \left. \liminf_{n \to +\infty} {n^{3/2} \over \rho_\star^{n-2k}}\mathbb P_{\tilde x_k, a_k} (
m_{n-2k}\circ \theta^k\geq 1, 1\leq \vert M_{0, n-2k  }\circ \theta^k\vert \leq K)\right]}
_
{\geq   c  \ h^{\theta_\star} (\tilde x_k, a_k)\  {\tilde h}^{\theta_\star} (\tilde x_k, 0) e^{\theta_\star(a_k-\ln K)}\ \ln K\ >0\ {\rm by}\  (\ref{localminimumliminf}) }
\\
& 
\geq c\  {  \ln K\over K^{\theta_\star}}  \ \   
\underbrace{\rho_\star^{-k} \mathbb E \left[{1\over  1+{\bf A} } \vert M_{0, k-1}\vert^{\theta_\star};  m_{k} \geq 1, h^{\theta_\star} (\tilde{\bf 1}\cdot M_{0, k-1}, \ln\vert M_{0, k-1}\vert
)\right]}_{\ \asymp \quad \widehat {\mathbb{E}}_{\tilde {\bf 1}, 0}^{\theta_\star}\left[{1\over  1+{\bf A} }\right] }
\end{align*}
On the other hand, since the function $\tilde h^{\theta_\star}$ associated with the process $(M'_{n, 1}\cdot x, b+ \ln\vert M'_{n, 1}  x\vert)_{n \geq 1}$  satisfies $\tilde h^{\theta_\star}(x, b)\asymp  1$ for $b$ in a compact set, it holds
 \begin{align*}
  \rho_\star^{-k} \mathbb E \left[{1\over  1+{\bf C} '};  {1\over K} \leq \vert  M'_{k-1, 1}\vert   \leq 1, {\bf M}'_{k-1} \leq 1 \right]
\geq  c \  (\ln K)  \ 
\widehat {\mathbb{E}}_{{\bf 1}, 0}^{\theta_\star}\left[{1\over  1+{\bf C}' }\right] 
\end{align*}
 Finally, by (\ref{hsqdkhfkla}),
   \begin{align}\label{liminffinal}
 \liminf_{n \to +\infty} {n^{3/2}\over \rho_\star^n} \mathbb E \left[{1\over {\bf A}+{\bf C}}; m_n \geq 1\right]
 \succeq (\ln K) \  \widehat{\mathbb{E}}_{\tilde {\bf 1}, 0}^{\theta_\star}\left[{1\over  1+{\bf A} }\right] \times
  \widehat {\mathbb{E}}_{{\bf 1}, 0}^{\theta_\star}\left[{1\over  1+{\bf C}' }\right]. 
 \end{align}
 
\noindent {\bf Second step: control of the term $\displaystyle    \mathbb E \left[{\bf B}; m_n \geq 1\right]$}
\begin{align*}
\mathbb E \left[{\bf B}; m_n \geq 1\right]
&= \sum_{\ell=k}^{n-k-1} \mathbb E\left[ \vert M_{0, \ell}\vert ^{-1} ;  m_n \geq 1\right]\\
&
\leq  \sum_{j \geq 0}  \sum_{\ell=k}^{n-k-1} 2^{-j} \mathbb P\left( \vert M_{0, \ell}\vert ^{-1}  \in [2^j, 2^{j +1}[;  m_n \geq 1\right)
\\
 & 
\leq  \sum_{j \geq 0}  \sum_{\ell=k}^{n-k-1} 2^{-j} \mathbb P( \vert M_{0, \ell}\vert ^{-1}  \in [2^j, 2^{j +1}[;  m_\ell\geq 1, 
\\
&\qquad \qquad \qquad \qquad \qquad \qquad  \vert M_{0, \ell+1}\vert \geq 1, \ldots,  \vert M_{0, n-1}\vert \geq 1)
\end{align*}
\begin{align*}
\qquad &\leq  \sum_{j \geq 0}  \sum_{\ell=k}^{n-k-1} 2^{-j} \mathbb P( \vert M_{0, \ell}\vert ^{-1}  \in [2^j, 2^{j+1}[;  m_\ell \geq 1, 
\\
&\qquad \qquad \qquad \qquad \qquad \qquad  \vert M_{\ell+1, \ell+1}\vert \geq  2^{-j}, \ldots,  \vert M_{\ell+1, n-1} \vert \geq 2^{-j} )
\\
& 
=  \sum_{j\geq 0}  \sum_{\ell=k}^{n-k-1} 2^{-j} 
\underbrace{\mathbb P( \vert M_{0, \ell}\vert ^{-1}  \in [2^j, 2^{j +1}[;  m_\ell \geq 1)}_{\asymp {\rho_\star^{\ell}\over \ell^{3/2}}}
 \underbrace{\mathbb P(m_{n-\ell-1} \geq 2^{-j})}_{\asymp j { \rho_\star^{n-\ell}\over (n-\ell)^{3/2}}}
 \\
 &
 \preceq {\rho_\star^n\over n^{3/2}}\ \left(\sum_{j \geq 0}j2^{-j}\right)\ \underbrace{\left(\sum_{\ell=k}^{n-k-1}{n^{3/2}\over \ell^{3/2} (n-\ell)^{3/2} }\right)}_{\preceq {1 \over \sqrt{k}}}
\end{align*}

 \noindent {\bf Last step: choosing $k$ large enough}
 
   By the above  $\displaystyle  \sup_{n \geq }{n^{3/2}\over \rho_\star^n}\mathbb E \left[{\bf B}; m_n \geq 1\right] \longrightarrow 0$ as $k \to +\infty$. It remains to check that the factors
$\displaystyle 
 \widehat {\mathbb{E}}_{\tilde {\bf 1}, 0}^{\theta_\star}\left[{1\over  1+{\bf A} }\right]$ and $\displaystyle   \widehat {\mathbb{E}}_{{\bf 1}, 0}^{\theta_\star}\left[{1\over  1+{\bf C}' }\right]$ 
%
 in the right hand side  of (\ref{liminffinal}) do not converge  to $0$ as $k \to +\infty$.
 
By (\ref{hsqdkhfkla}) and by convexity of the function $x \mapsto {1\over x}$ on $\mathbb R^{*+}$, it holds
 \[  \widehat {\mathbb{E}}_{\tilde {\bf 1}, 0}^{\theta_\star}\left[{1\over 1+{\bf A}}\right]  \geq   {1\over \widehat {\mathbb{E}}_{\tilde {\bf 1}, 0}^{\theta_\star} [1+{\bf A}]}\geq {1\over \widehat {\mathbb E}_{\tilde {\bf 1}, 0}^{\theta_\star} \displaystyle  \left[1+\sum_{\ell=0}^{+\infty} \vert \tilde{\bf 1}M_{0, \ell}\vert^{-1}\right]}.
 \]
 and it suffices to check that $\displaystyle \widehat {\mathbb E}_{\tilde {\bf 1}, 0}^{\theta_\star} \displaystyle  \left[\sum_{\ell=0}^{+\infty} \vert \tilde{\bf 1}M_{0, \ell}\vert^{-1}\right]<+\infty$. Indeed, for any $\ell \geq 0$,
\begin{align*}
 \widehat {\mathbb E}_{\tilde {\bf 1}, 0}^{\theta_\star}\left [\vert \tilde{\bf 1}M_{0, \ell}\vert^{-1}\right]
 &\preceq 
 {\mathbb E} _{\tilde {\bf 1}, 0}^{\theta_\star}\left [\vert  M_{0, \ell}\vert^{-1}h^{\theta_\star}(X_{\ell+1}, S_{\ell+1}), \tau > \ell +1 \right]\\
 &\preceq 
 {\mathbb E} ^{\theta_\star}\left [\vert  M_{0, \ell}\vert^{-1}(1+\ln \vert   M_{0, \ell}\vert), m_{\ell +1} \geq 1\right]\\
 &\leq  \sum_{j=0}^{+\infty}(1+j) e^{- j} \  \underbrace{\mathbb P^{\theta_\star}( e^j\leq \vert M_{0, \ell}\vert <e^{j+1}, m_{\ell +1} \geq 1)}_{=\ \mathbb P^{\theta_\star}_{\tilde{\bf 1}, 0}(j\leq S_{\ell +1} <j+1, \tau >\ell +1)}\\
  &\leq {1\over \ell^{3/2}}\sum_{j=0}^{+\infty} (1+j)^2e^{- j}  \quad {\rm by \  inequality} \ (\ref{3/2upper}),
\end{align*}
which yields  the expected result.
The same argument works for the factor $\displaystyle   \widehat {\mathbb{E}}_{{\bf 1}, 0}^{\theta_\star}\left[{1\over  1+{\bf C}' }\right]$ 
by considering the random variables $ \tilde S_{\ell}({\bf 1}, 0)$ and $\tilde \tau_{{\bf 1}, 0}$ instead of  $S_{\ell}(\tilde{\bf 1}, 0)$ and 
$ \tau_{\tilde{\bf 1}, 0}$.

\end{document}